\documentclass[12pt]{amsart}
\usepackage{amsmath, amssymb,graphicx,caption, subcaption, mathtools, enumitem}
\usepackage[pdftex,hyperindex=true]{hyperref}
\usepackage{url}
\usepackage{etoolbox}
\patchcmd{\section}{\scshape}{\bfseries}{}{}
\patchcmd{\subsection}{\bfseries}{\itshape}{}{}

\makeatletter
\def\@seccntformat#1{%
  \protect\textup{\protect\@secnumfont
    \ifnum\pdfstrcmp{section}{#1}=0 \bfseries\fi% section # in \bfseries
    \ifnum\pdfstrcmp{subsection}{#1}=0 \itshape \fi% subsection # in \itshape
    \csname the#1\endcsname
    \protect\@secnumpunct
  }%
}  
\makeatother

\usepackage[utf8]{inputenc}
\usepackage[english]{babel}

\theoremstyle{plain}

\newtheorem*{conjecture}{Conjecture}

\theoremstyle{definition}

\newtheoremstyle{note}% <name>
{3pt}% <Space above>
{3pt}% <Space below>
{\itshape}% <Body font>
{}% <Indent amount>
{\itshape}% <Theorem head font>
{:}% <Punctuation after theorem head>
{.5em}% <Space after theorem headi>
{}% <Theorem head spec (can be left empty, meaning `normal')>

\theoremstyle{note}

\theoremstyle{plain} % just in case the style had changed
\newcommand{\thistheoremname}{}
\newtheorem*{genericthm}{\thistheoremname}

\numberwithin{equation}{section}

\newcommand{\acknowledge}{\subsection*{Acknowledgements}}

\def\Dbar{\leavevmode\lower.6ex\hbox to 0pt{\hskip-.23ex \accent"16\hss}D}

\begin{document}

\title{Root Systems and the Atiyah-Sutcliffe Problem}
\author{J. Malkoun}
\address{Department of Mathematics and Statistics\\
Faculty of Natural and Applied Sciences\\
Notre Dame University-Louaize, Zouk Mosbeh\\
P.O.Box: 72, Zouk Mikael,
Lebanon}
\email{joseph.malkoun@ndu.edu.lb}

\date{Received: date / Accepted: date}
\maketitle

\begin{abstract} In this short note, we show that the Atiyah-Sutcliffe conjectures for $n = 2m$, 
related to the unitary groups $U(2m)$, imply the author's analogous conjectures, which are associated with the symplectic 
groups $Sp(m)$. The proof is based on the simple fact that the root system of $U(2m)$ dominates that of $Sp(m)$.
\end{abstract}

\maketitle

\section{Introduction} \label{intro}

In \cite{BR1997}, M.V. Berry and J.M. Robbins considered a new argument for the spin-statistics theorem in Quantum Mechanics. Their motivation 
was that the standard arguments for this theorem usually involves notions belonging to Quantum Field Theory, while the spin-statistics theorem itself 
is a quantum-mechanical statement. In \cite{BR1997}, the sign
\[ (-1)^{\frac{s(s-1)}{2}} \]
which the wave function of a collection of $n$ identical spin $s$ particles picks up upon interchanging two of the particles, arises as a 
geometric Berry phase factor. However the discussion in \cite{BR1997} was mostly for two and three identical particles. To extend 
their work to $n$ identical particles, the authors needed the existence of a certain map. This led to the formulation of the Berry-Robbins problem, 
which we will now state.

Define $C_n(\mathbb{R}^3)$ to be the configuration space of $n$ distinct points in $\mathbb{R}^3$. Denote by $U(n)/T^n$ the symmetric space 
obtained as the quotient of the unitary group $U(n)$ by a maximal torus $T^n$, which we will take for simplicity to consist of all diagonal 
$n$-by-$n$ matrices with entries in $U(1)$.

The symmetric group $\Sigma_n$ on $n$ letters acts on $C_n(\mathbb{R}^3)$ by permuting the $n$ (distinct) points in $\mathbb{R}^3$ in every 
configuration, and it acts on $U(n)$ by permuting the columns of every $n$-by-$n$ unitary matrix. The latter action descends to an action of 
$\Sigma_n$ on $U(n)/T^n$.
\\
\smallskip
\\
\emph{Berry-Robbins Problem:} Does there exist, for each $n \geq 2$, a continuous map 
\[ f_n: C_n(\mathbb{R}^3) \to U(n)/T^n \]
which is $\Sigma_n$-equivariant (for the action of $\Sigma_n$ defined in the previous paragraph)?
\\
\smallskip
\\
The Berry-Robbins problem was solved positively by M.F. Atiyah in \cite{Ati-2000}. Nevertheless, M.F. Atiyah was not satisfied with that solution, 
and proposed another candidate map in \cite{Ati-2001}, which had more desirable features than the maps in \cite{Ati-2000}, being in particular 
\emph{smooth}, and also $SO(3)$-equivariant (in addition to being $\Sigma_n$-equivariant of course). However, for them to be genuine solutions of the 
Berry-Robbins problem, a linear independence statement had to hold. Let us refer to the latter as the linear independence conjecture, denoted simply 
by LIC.

The LIC was further strengthened by M.F. Atiyah and P.M. Sutcliffe to two successively stronger statements in \cite{Ati-Sut-2002}, using a notion of 
normalized determinant function, referred to as the Atiyah-Sutcliffe determinant. The LIC, together with these two stronger conjectures, are referred to 
collectively as the Atiyah-Sutcliffe conjectures.

In his Edinburgh Lectures on Geometry, Analysis and Physics, M.F. Atiyah asked whether there exist similar maps corresponding to other Lie groups. 
As a partial answer to this question, the author of this article found in \cite{Malkoun2014} a variant of the Atiyah-Sutcliffe maps for the symplectic groups 
$Sp(m)$, instead of the unitary groups $U(n)$. Moreover, conjectures similar to the Atiyah-Sutcliffe conjectures were formulated in \cite{Malkoun2014}. Let 
us refer to this collection of conjectures as the symplectic conjectures.

The aim of the present note is to show that the Atiyah-Sutcliffe conjectures corresponding to $U(2m)$ imply the symplectic conjectures corresponding to 
$Sp(m)$. The key observation of the proof is that the root system $U(2m)$ dominates that of $Sp(m)$.

\section{Review of the Atiyah-Sutcliffe Problem on Configurations of Points}

Let $\mathbf{x} \in C_n(\mathbb{R}^3)$, and write $\mathbf{x} = (\mathbf{x}_1,\ldots,\mathbf{x}_n)$, where $\mathbf{x}_a \in \mathbb{R}^3$ for 
every $a$, $1 \leq a \leq n$. Consider the Hopf map
\[ h: \mathbb{C}^2 \setminus \{0\} \to \mathbb{R}^3 \setminus \{0\} \]
defined by
\begin{equation} \label{hopf} h(u,v) = (2u \bar{v}, |u|^2-|v|^2) \in \mathbb{C} \times \mathbb{R} \setminus \{ (0,0) \} \end{equation}
where $\mathbb{C} \times \mathbb{R}$ is identified with  $\mathbb{R}^3$ by identifying $\mathbb{C}$ with $\mathbb{R}^2$ in the natural way (map $\zeta \in \mathbb{C}$ 
to $(x,y) \in \mathbb{R}^2$, where $x$, resp. $y$, is the real, resp. imaginary, part of $\zeta$).

For each pair of indices $a$, $b$, where $1 \leq a,b \leq n$, $a \neq b$, we form the direction vector
\[ \mathbf{x}_{ab} = \mathbf{x}_b - \mathbf{x}_a \]
We then choose a Hopf lift $\mathbf{u}_{ab} = (u_{ab},v_{ab}) \in \mathbb{C}^2 \setminus \{0\}$ of $\mathbf{x}_{ab}$, so that
\[ h(\mathbf{u}_{ab}) = \mathbf{x}_{ab} \]
A Hopf lift is of course only defined up to multiplication by a scalar in $U(1)$ ($h(\lambda \mathbf{u}) = h(\mathbf{u})$ for any $\lambda \in U(1)$).
Once the Hopf lifts are all made, we then form for each pair of indices $a$, $b$ ($1 \leq a,b \leq n$, $a \neq b$), the following polynomials
\[ p_{ab}(t) = u_{ab} t - v_{ab} \]
where $t$ is a complex variable. Finally, we form the following polynomials
\[ p_a(t) = \prod_{\substack{b \neq a \\  1 \leq b \leq n}} p_{ab}(t) \qquad \text{($1 \leq a \leq n$)}\]

We can now state the LIC.

\begin{conjecture}[LIC] The polynomials $p_a$, for $1 \leq a \leq n$, are linearly independent over $\mathbb{C}$. \end{conjecture}

We denote by $(p_1,\ldots,p_n)$ the $n$-by-$n$ matrix having as $j$'th column the coefficients of the polynomial $p_j$, with the convention that they are ordered according 
to increasing powers of $t$. Similarly, we denote by $(p_{ab}, p_{ba})$ the $2$-by-$2$ matrix having as first (resp. second) column the coefficients of $p_{ab}$ (resp. $p_{ba}$), following 
a similar ordering convention on the coefficients of each polynomial. The Atiyah-Sutcliffe (normalized) determinant map
\[ D: C_n(\mathbb{R}^3) \to \mathbb{C} \]
can now be defined by
\begin{equation} D(\mathbf{x}) = \frac{\operatorname{det}(p_1,\ldots,p_n)}{\prod_{1\leq a<b \leq n} \operatorname{det}(p_{ab},p_{ba})} \label{ASdef} \end{equation}

It is clear from the definition \eqref{ASdef} of the Atiyah-Sutcliffe determinant $D$ that the LIC is equivalent to the non-vanishing of the map $D$, i.e. $D(\mathbf{x}) \neq 0$ 
for every $\mathbf{x} \in C_n(\mathbb{R}^3)$. In \cite{Ati-Sut-2002}, the authors actually conjectured a stronger statement (known as conjecture 2 in that article), namely the following.

\begin{conjecture}[Conjecture 2] For every $\mathbf{x} \in C_n(\mathbb{R}^3)$, $|D(\mathbf{x})| \geq 1$. \end{conjecture}

The authors of \cite{Ati-Sut-2002} also formulated an even stronger conjecture, referred to as Conjecture 3 in \cite{Ati-Sut-2002}, but we shall refer the interested reader to that 
article. The LIC, together with conjecture 2 and 3 are referred to collectively as the Atiyah-Sutcliffe conjectures.

\section{The Symplectic Conjectures}

We present the author's variant of the Atiyah-Sutcliffe problem on configurations of points, related to the symplectic groups $Sp(m)$, instead of the unitary groups. Define
\begin{align*} \mathcal{C}_m(\mathbb{R}^3) = \{ \mathbf{x} = (\mathbf{x}_1,\ldots,\mathbf{x}_m) \in (\mathbb{R}^3)^m;& \mathbf{x}_a \neq \mathbf{0} \text{ for $1\leq a \leq m$ and } \\ 
&\mathbf{x}_a \pm \mathbf{x}_b \neq \mathbf{0} \text{ for $1 \leq a<b \leq m$} \} \end{align*}
We also define the set
\[ I = \{1,\ldots,m\} \cup \{\bar{1},\ldots,\bar{m}\} \]
Let us define a total order on $I$ by
\begin{equation} 1 < \bar{1} < \cdots < m < \bar{m} \label{order} \end{equation}
For each $\alpha \in I$, we define
\[ \mathbf{x}_\alpha = \left\{ \begin{array}{cc} \mathbf{x}_a, &\text{if $\alpha = a$, for $1 \leq a \leq m$} \\
                                                   -\mathbf{x}_a, &\text{if $\alpha = \bar{a}$, for $1 \leq a \leq m$} \end{array} \right. \]
                                                   
We then choose, for each $\alpha,\beta \in I$, $\alpha \neq \beta$, a Hopf lift 
$\mathbf{u}_{\alpha\beta} \in \mathbb{C}^2 \setminus \{0\}$ of $\mathbf{x}_{\beta}-\mathbf{x}_{\alpha}$, i.e.
\[ h(\mathbf{u}_{\alpha\beta}) = \mathbf{x}_{\beta} - \mathbf{x}_{\alpha} \]
We remark that $\mathbf{x}_{\beta} - \mathbf{x}_{\alpha} \in \mathbb{R}^3 \setminus \{ 0 \}$ since $\mathbf{x} \in \mathcal{C}_m(\mathbb{R}^3)$.
Once the Hopf lifts are chosen, we then form, for each pair of indices $\alpha, \beta \in I$, with $\alpha \neq \beta$, the polynomials
\[ p_{\alpha\beta}(t) = u_{\alpha\beta} t - v_{\alpha\beta} \]
where $t$ is a complex variable.
For each $\alpha \in I$, we form the polynomial
\begin{equation} p_{\alpha}(t) = \prod_{\beta \in I \setminus \{\alpha\}} p_{\alpha \beta}(t) \end{equation}

\begin{conjecture}[Symplectic LIC] The polynomials $p_{\alpha}$, for $\alpha \in I$, are linearly independent over $\mathbb{C}$. \end{conjecture}

We now denote by $(p_1,p_{\bar{1}},\ldots,p_m,p_{\bar{m}})$ the $2m$-by-$2m$ matrix having the coefficients of $p_1$ as first column, 
and so on, where the coefficients of each polynomial $p_\alpha$ are ordered by increasing powers of $t$.

We now define the symplectic analogue $D_S: \mathcal{C}_m(\mathbb{R}^3) \to \mathbb{C}$ of the Atiyah-Sutcliffe determinant function by
\[ D_S(\mathbf{x}) = \frac{\operatorname{det}(p_1,p_{\bar{1}},\ldots,p_m,p_{\bar{m}})}{\prod_{1 \leq \alpha < \beta \leq \bar{m}} \operatorname{det}(p_{\alpha \beta},p_{\beta \alpha})} \]
where the order relation $<$ is the one defined in \eqref{order}. 
  
The analogue of the Atiyah-Sutcliffe Conjecture 2 can now be formulated.

\begin{conjecture}[Symplectic Conjecture 2] For every $\mathbf{x} \in \mathcal{C}_m(\mathbb{R}^3)$, 
\[ |D_S(\mathbf{x})| \geq 1 .\] \end{conjecture}

While an analogue of the Atiyah-Sutcliffe Conjecture 3 can also be formulated for the symplectic variant of the original problem, we shall refrain from 
discussing it, for brevity (the interesting reader can fill in the details if they wish to do so).

\section{Proof of the Main Theorem}

Consider the following linear mapping
\[ g: \mathbb{R}^m \to \mathbb{R}^{2m} \]
mapping $e_a$ ($1 \leq a \leq m$) to $v_{2a-1} - v_{2a}$, where the $e_a$ (resp. the $v_{\alpha}$) are the 
canonical basis vectors of $\mathbb{R}^m$ (resp. $\mathbb{R}^{2m}$). Let us denote by $x_a$, for $1 \leq a \leq m$ 
(resp. $y_{\alpha}$, for $1 \leq \alpha \leq 2m$), the linear forms on $\mathbb{R}^m$ (resp. $\mathbb{R}^{2m}$) which form 
a dual basis of the $e_a$ (resp. the $v_{\alpha}$).

Denote by $A$ the following subset of the dual space of $\mathbb{R}^{2m}$
\[ A = \{ y_{\alpha}-y_{\beta}; 1 \leq \alpha, \beta \leq 2m \text{ and } \alpha \neq \beta \} \]
and by $C$ the following subset of the dual space of $\mathbb{R}^m$
\[ C = \{ \pm x_a \pm x_b; 1 \leq a < b \leq m \} \cup \{ \pm 2 x_a; 1 \leq a \leq m \}. \]
The sets $A$ and $C$ are actually the root systems of $U(2m)$ and $Sp(m)$ respectively. Then it can be easily checked that
\begin{equation} g^*(A) = C .\end{equation}
The linear map $g$ can be modified to a map 
\[ \hat{g}: \mathcal{C}_m(\mathbb{R}^3) \to C_{2m}(\mathbb{R}^3) \]
mapping $(\mathbf{x}_1,\ldots,\mathbf{x}_m) \in \mathcal{C}_m(\mathbb{R}^3)$ to
\[ (\mathbf{x}_1,-\mathbf{x}_1,\ldots,\mathbf{x}_m,-\mathbf{x}_m) \in C_{2m}(\mathbb{R}^3). \]
It can now be easily checked that the $2m$ polynomials $p_{\alpha}$, for $1 \leq \alpha \leq 2m$, associated to a 
configuration $\mathbf{x} \in \mathcal{C}_m(\mathbb{R}^3)$ via the symplectic maps, are the same (each up to a 
scalar factor) as the polynomials associated to the configuration $\hat{g}(\mathbf{x}) \in C_{2m}(\mathbb{R}^3)$ via 
the Atiyah-Sutcliffe maps. It then follows that
\begin{equation} D_S(\mathbf{x}) = D(\hat{g}(\mathbf{x})) .\end{equation}
Hence the symplectic conjectures for $Sp(m)$ are implied by the Atiyah-Sutcliffe conjectures for $U(2m)$, as claimed.

\acknowledge{The author would like to thank Sir Michael Atiyah for listening to his ideas. He would also like to thank John Stembridge for a nice 
discussion on the notion of folding Dynkin diagrams, and Hassan Azad for mentioning to him the notion of twisting.}

\smallskip

\vspace{5mm}

\def\Dbar{\leavevmode\lower.6ex\hbox to 0pt{\hskip-.23ex \accent"16\hss}D}
\providecommand{\bysame}{\leavevmode\hbox to3em{\hrulefill}\thinspace}
\providecommand{\MR}{\relax\ifhmode\unskip\space\fi MR }
\providecommand{\MRhref}[2]{%
  \href{http://www.ams.org/mathscinet-getitem?mr=#1}{#2}
}
\providecommand{\href}[2]{#2}

\end{document}